\documentclass[12pt]{article}
\usepackage{amsmath,amsfonts,amssymb,mathptmx}
\usepackage{mordell,amscd,theorem}
\usepackage{wrapfig}
\usepackage[all]{xy}

\newtheorem{theorem}{Theorem}[section] 

{\theorembodyfont{\rmfamily}   }
\newtheorem{lemma}[theorem]{Lemma}

{\theorembodyfont{\rmfamily}   }
{\theorembodyfont{\rmfamily}   }
{\theorembodyfont{\rmfamily}   }

\begin{document}

\begin{center}\huge
Note on the rank of quadratic twists \\
of Mordell equations\\
\large
\VS{1.5\baselineskip}Sungkon Chang\\
\normalsize
\texttt{changsun@mail.armstrong.edu}
\end{center}

\begin{abstract} 
Let $E$ be the elliptic curve given by a Mordell equation $y^2=x^3-A$ where $A \in \zz$.  
Michael Stoll found a precise formula for the size of a Selmer group of  $E$ for certain values of $A$.
For $D\in\zz$, let $E_D$ denote the \textit{quadratic twist} $D\, y^2 =x^3-A$.
We use Stoll's formula to show that for a positive square-free integer  $A \equiv 1$ or $25 \mod 36$ and for a non-negative integer $k$, we can compute a lower bound for  the proportion of square-free integers $D$ up to $X$ such that 
$\rank E_D(\ratn) \le 2k$.
We also compute an upper bound for a certain average rank of quadratic twists of $E$.
\end{abstract}

\section{Introduction}\label{sec:intro}

Let $E/\ratn$ be the elliptic curve $ y^2=x^3-A$ where $A$ is a nonzero integer.  Then,  let us denote by  $E_D$  the \textit{quadratic twist} $y^2=x^3-AD^3$ for each nonzero square-free integer $D$.
\paragraph{}
 For a non-negative integer $k$, let 
$$ \delta_k:=\frac{ 3^{k+1} - 2 }{  3^{k+1} - 1 },$$
and let $T(X)$ denote the set of all positive square-free integers less than $X$.
In this paper, we shall prove the following result:
\par
\begin{copythm}{Theorem \ref{thm:MainTheorem1}}
Let $A$ be a  positive square-free  integer such that $A\equiv 1$ or $25 \mod 36$. Then, 
for a non-negative integer $k$, 
\begin{equation}\label{eq:thm1-1}
\lim\inf_X\ \frac{ \Sharp\ \set{ D\in T(X) :  \rank E_D(\ratn)\le 2k } }{  \Sharp\ T(X) } \ge 
         \frac{\delta_k}{8} \cdot \prod_{p \mid A} \frac{p}{(p-1)(p+1)}.
\end{equation}
In particular, 
\begin{equation}\label{eq:thm1-2}
\lim\inf_X\ \frac{ \Sharp\ \set{ D\in T(X) :  E_D(\ratn)=\set{O} } }{  \Sharp\ T(X) } \ge 
         \frac{1}{16} \cdot \prod_{p \mid A} \frac{p}{(p-1)(p+1)}.
\end{equation} 
\end{copythm}
\par
One of the earliest (known) examples of elliptic curves $E/\ratn$ with a positive proportion of square-free integers $D$ such that $\rank E_D(\ratn)=0$ is the elliptic curve given by $y^2=x^3-x$, proved by the work of Heath-Brown in \cite{heath-brown:1994}, 1994.
Note that a similar result was already available in the late eighties.  
There are two  results known in 1988 which together simply imply that the elliptic curve $E : y^2=x^3-1$ has positive proportion of quadratic twists of rank $0$.  
In 1985, Frey proved in \cite[Proposition, p.237]{frey:1985} that if $D$ is a square-free integer such that $D\equiv 1 \mod 4$, then 
\begin{center}
$\Sharp\Cl(\ratn(\st{-D}))[3]=1$ if and only if $\Sel{3}(E_D,\ratn)=\set{0}$
\end{center}
where $\Sel{3}(E_D,\ratn)$ is \textit{the $3$-Selmer group of $E_D/\ratn$.}
In 1988, Nakagawa and Horie 
 proved in \cite{nakagawa:1988} Theorem \ref{thm:davenport} stated in this paper, which is a refined result   of the famous theorem of Davenport and Heilbronn.  Their theorem implies that there is a positive proportion of square-free positive integers $D$ such that $D\equiv 1 \mod 4$ and  $\Sharp\Cl(\ratn(\st{-D}))[3]=1$. 
Therefore, it follows that there is a positive proportion of (positive) square-free integers $D$ such that 
\begin{center}
$\Sel{3}(E_D,\ratn)=\set{0}$ and, hence, $\rank E_D(\ratn)=0$.
\end{center}

\paragraph{}
Let us introduce our second result. 
For two positive integers $m$ and $N$, let us denote by $\Ntwop(X,m,N)$ the set of positive \textit{fundamental discriminants} $\Delta < X$ such that $\Delta \equiv m \mod N$.

\begin{copythm}{Theorem \ref{thm:MainTheorem2}}
Let $E/\ratn$ be an elliptic curve given by $y^2=x^3-A$ for some $A \in \zz$ such that 
 $A \equiv 1$ or $25 \mod 36$ is a square-free integer.
\par
If $A>0$, then 
$$
\lim\sup_{X\to \infty} \frac{ \sum_{D \in \Ntwop(X,1,12A)} \rank( E_D(\ratn)) }{ \Sharp\ \Ntwop(X,1,12A) } 
%\le
%\lim_{X\to \infty} \frac{ \sum_{D \in \Ntwop(X,1,36A)} \dimp3 \Sel{\lambda}(E_D/K) }{ \Sharp\ \Ntwop(X,1,36A) }
\le 1.$$
If $A<0$, then 
$$
\lim\sup_{X\to \infty} \frac{ \sum_{D \in \Ntwop(X,1,12\abs{A})} \rank( E_D(\ratn)) }{ \Sharp\ \Ntwop(X,1,12\abs{A}) } 
%\le 
%\lim_{X\to \infty} \frac{ \sum_{D \in \Ntwop(X,1,36A)} \dimp3 \Sel{\lambda}(E_D/K) }{ \Sharp\ \Ntwop(X,1,36A) }
                \le \frac43.$$
\end{copythm}
\par
Recall that $T(X)$ denotes the set of all positive square-free integers less than $X$.
Assuming the Birch and Swinnerton-Dyer Conjecture, (the Modularity Conjecture), and a form of the Riemann Hypothesis, 
Goldfeld proved in \cite{goldfeld:1979} that 
\begin{equation}\label{eq:average}
\underset{X\to\infty}{\lim\sup}\ 
        \frac{ \sum_{ \abs{D}\in T(X) } \rank E_D(\ratn) }{ 2\cdot \Sharp\ T(X) } \le 3.25.
        \end{equation}
In \cite{heath:2004}, this  upper bound is reduced to $1.5$ by Heath-Brown.
In \cite{goldfeld:1979}, Goldfeld conjectured that the average in (\ref{eq:average}) should be $1/2$, which is known as Goldfeld's conjecture.
In \cite{heath-brown:1994}, Heath-Brown (unconditionally) computes an upper bound for the average rank of quadratic twists  $y^2=x^3-D^2\,x$ over odd integers $D$, and  Gang Yu in \cite{yu:2003} computes a certain average rank of quadratic twists of infinitely many elliptic curves with rational $2$-torsion points. At this moment of writing, these two results were the only unconditional results, known to the author, on the average rank of quadratic twists of an elliptic curve.
\paragraph{Acknowledgements}
I would like to thank Professors Roger Heath-Brown,  
Dino Lorenzini, and Michael Stoll for helpful comments.
I also wish to thank Ken Ono, Kevin James, Su-ion Ih, and Charles Pooh for their encouragement on this work.

\subsection{Stoll's Formula}\label{sec:formula}
Let $\zeta\in\Qbar$ be a primitive third root of unity, and let $\lambda:=1-\zeta$. Let $K$ denote the cyclotomic field extension $\ratn(\zeta)$.
Let $E/K$ be an elliptic curve given by $y^2=x^3-A$ where $A \in \zz$.  We denote simply by $\zeta$ the endomorphism on $E$ given by $(x,y) \mapsto (\zeta x ,y)$ which is defined over $K$.  Let us denote the endomorphism $1-\zeta$ on $E$ simply by $\lambda$, and let $E[\lambda]$ denote  the kernel of $\lambda$.
The endomorphism $\lambda$ induces the Kummer sequence
\begin{equation}\label{eq:kummer-seq}
 0 \To E[\lambda] \To E \overset{\lambda}{\To} E \To 0.
 \end{equation}
Let $\MK$ denote the set of all places of $K$, and let $K_v$ denote the completion of $K$ with respect to  a place $v\in \MK$.  
Note that (\ref{eq:kummer-seq}) induces the following injective homomorphism into the first cohomology of the $\Gal(\overline F/F)$-module $E[\lambda]$ where $F$ is the number field $K$ or a completion $K_v$:
\begin{equation}\label{eq:fund-seq}
\delta_F : E(F)/\lambda E(F) \To \HonE(F,E[\lambda]).
\end{equation}
 When  $F=K_v$ for some $v \in \MK$, let us denote $\delta_F$ also by $\delta_v$.
\par
Note that for each $v \in \MK$, there is the restriction map $\res_v : \HonE(K,E[\lambda]) \to \HonE(K_v,E[\lambda])$ (see \cite[Chapter X, Sec 4]{silverman:1986}).
\textit{The $\lambda$-Selmer group of $E/K$} is 
\begin{equation}\label{eq:def-Sel}
\Sel{\lambda}(E,K):=\set{ \xi \in \HonE(K, E[\lambda]) : \res_v(\xi) \in \Img \delta_v \text{ for all } v \in \MK},
\end{equation}
and it contains the image of $E(K) / \lambda E(K)$.

\begin{theorem}\label{prop:stoll}\ThmAuthor{Stoll, \cite[Corollary 2.1]{stoll:1998}}
Let $A$ be a rational integer. Let $E/K$ be the elliptic curve given by $y^2 = x^3 - A$.
\par
Suppose that the following conditions\footnote{In \cite{stoll:2002}, which is a sequel of \cite{stoll:1998}, he improved the conditions so that more possibilities of values of $A$ can be considered for Theorem \ref{thm:MainTheorem1}.}
are satisfied:
\begin{enumerate}
\item\label{prop:cond1}  $-A\equiv 2 \mod 3$.
\item\label{prop:cond2} For all places $v \ne \lambda$ of $K$ of bad reduction for $E/K$, the integer $-A$ is non-square in $K_v^*$, e.g., $-A$ is square-free, and $-A\equiv 3 \mod 4$.
\end{enumerate}
Then, $\dimp3\sellam(E,K)=$
\begin{alignat}{2}
1 + &2 \dimp3 \Cl(\ratn(\sqrt{-A}))[3] &&\quad\text{ if } -A \equiv 2,8 \mod 9 \text{ and } A<0,\notag\\
       & 2 \dimp3 \Cl(\ratn(\sqrt{-A}))[3] &&\quad\text{ if } -A \equiv 2,8 \mod 9 \text{ and } A>0,\notag\\
       & 2 \dimp3 \Cl(\ratn(\sqrt{3A}))[3] &&\quad\text{ if } -A \equiv 5 \mod 9 \text{ and } A<0,\notag\\     
1 + & 2 \dimp3 \Cl(\ratn(\sqrt{3A}))[3] &&\quad\text{ if } -A \equiv 5 \mod 9 \text{ and } A>0.\notag
\end{alignat}
In particular, these numbers give a bound on $\rank E(\ratn)$.
\end{theorem}

\begin{lemma}\label{lem:condition}
If $n$ is a nonzero integer and $p$ is a prime number $> 3$ such that $n \equiv 3 \mod 4$ and $\ord_p(n) \equiv 1 \mod 2$, then
 $n \not\in (\Kp^*)^2$ where $\Kp$ is the completion at any prime ideal $\primep$ of $\OK$ lying over $p$ or $2$.
\par
Let $A$ be a square-free integer $\equiv 1 \mod 12$.
If $D$ is a square-free integer coprime to $A$ such that $D \equiv 1 \mod 12$,  then the elliptic curve: $y^2=x^3-A D^3$ satisfies Conditions (\ref{prop:cond1}) and (\ref{prop:cond2}) in Theorem \ref{prop:stoll}.
\end{lemma}
\proof
The proof is left to the reader.\hfill\qed
\par
In this paper, we will focus on quadratic twists of the elliptic curve given by a Mordell equation, but the reader might have noticed from the formula in Theorem \ref{prop:stoll} that if $A$ is replaced with $AD^2$ for an integer $D$ such that $D \equiv 1 \mod 9$, and such that the elliptic curve $E^D : y^2=x^3-AD^2$ satisfies the conditions required for the formula, then the size of the Selmer group of $E^D/K$ equals that of the Selmer group of $E/K$. Since $y^2 = x^3 - AD^2$ forms a family of cubic twists, we can use the formula to obtain the following result on the distribution of  Mordell-Weil rank of cubic twists of $E$.
If $A$ is a positive square-free integer such that $A\equiv 1$ or $25 \mod 36$ and $\dimp3 \Cl(\ratn(\sqrt{-A}))[3]=0$, then there is a positive real number $\ep <1$ such that 
\begin{equation}\label{eq:cubic}
 \Sharp\, \set{ 0< D < X : D\text{ cube-free,}\  \rank E^D(\ratn) = 0 } \gg \frac{ X }{ (\log X)^\ep }.
 \end{equation}
 To compute the lower bound in (\ref{eq:cubic}), we construct a set of prime numbers with positive Dirichlet density, and show that whenever $D$ is a positive integer divisible only by prime numbers contained in this set, the Mordell-Weil rank of  $E^D$ is $0$.
 This observation is generalized for \textit{superelliptic curves over global fields} in \cite{chang-thesis:2004}.
  To my knowledge, the only known example of an elliptic curve with infinitely many cubic twists of Mordell-Weil rank $0$ is $x^3+y^3=D$ proved by D.~Lieman \cite{lieman:1994}.
\par
By using the proof of the main theorems in \cite{chang-thesis:2004}, we can find an explicit congruence condition on prime numbers such that if $D$ is supported by these prime numbers, then the cubic/quadratic twist $E_D$ has Mordell-Weil rank $0$.  This result is proved without using the modularity of elliptic curves.  Hence, it is an interesting project to show that the Hasse-Weil $L$-functions of these twists really do not vanish at $s=1$.
By showing this, we can establish the fact that there are infinitely many cubic/quadratic twists for which a part of Birch and Swinnerton-Dyer Conjecture is true, namely, Mordell-Weil rank being $0$ for these examples implies (with no use of the modularity)  the non-vanishing of the $L$-functions at $s=1$.

\subsection{The refined result of  Davenport-Heilbronn}\label{sec:davenport}
In \cite{nakagawa:1988}, Nakagawa and Horie  proved a refined result of Davenport and Heilbronn. Let $N$ and $m$ be positive integers.
Let $\Ntwo(X,m,N)$ be the set of  fundamental discriminants $\Delta$ such that $-X < 
\Delta <0$, and $\Delta \equiv m \mod N$, and let $\Ntwop(X,m,N)$ be the set of  fundamental discriminants $\Delta$ such that $0 < 
\Delta <X$, and $\Delta \equiv m \mod N$.
Let $h_3(\Delta)$ denote $\Sharp\ \Cl(F)[3]$ where $F$ is the quadratic extension of $\ratn$ with discriminant $\Delta$.
\par
Let us describe the property for $N$ and $m$, which  we require for Theorem \ref{thm:davenport}.
\begin{quote}\textbf{Condition (**)}\ \
If an odd prime number $p$ is a common divisor of $m$ and $N$, then $p^2 \mid N$ and $p^2 \nmid m$. Further if $N$ is even, then  $4 \mid N$ and $m \equiv 1 \mod 4$,  or  $16 \mid N$ and $m \equiv 8 $ or $12 \mod 16$.
\end{quote}
\begin{theorem}\label{thm:davenport}\ThmAuthor{Nakagawa-Horie, \cite{nakagawa:1988}}
Let $N$ and $m$ be positive integers satisfying Condition (**). 
Then, 
\begin{gather}
\lim_{X\to\infty}\frac{1}{\Sharp\ \Ntwop(X,m,N)}\ \sum_{\Delta \in \Ntwop(X,m,N)} h_3(\Delta) = \frac{4}{3},\\
\lim_{X\to\infty}\frac{1}{\Sharp\ \Ntwo(X,m,N)}\ \sum_{\Delta \in \Ntwo(X,m,N)} h_3(\Delta) = 2.
\end{gather}
\end{theorem}

\section{Proof of Theorem \ref{thm:MainTheorem1}}

 Let $S$ be a subset of $\zz$, and for a positive integer $x$, let $S(x)$ denote the set of integers $n$ contained in $S$ such  that $\abs{n}<x$.
Let $\nat$ be the set of positive integers. 
Let $h$ be a set-theoretic function: $\nat \to \nat$ such that the images of $h$ are powers of $3$.  
For a non-negative integer $k$, let 
\begin{equation}\label{eq:Skx}
S_k(x) := \set{ a \in S(x) : h(a) \le 3^k } \text{ and } \delta_k(x):= \frac{\SSkx}{\SSx}.
\end{equation}

\begin{lemma}\label{lem:liminf}
If $\limx \frac{1}{\SSx} \sum_{a \in \Sx} h(a) = B$ for some positive real number $B$, then
for a non-negative integer $k$,
$$ \underset{x\to\infty}{\lim \inf}\ \delk \ge \frac{ 3^{k+1} - B }{ 3^{k+1} - 1 }.$$
\end{lemma}

\proof
The proof is left to the reader.\hfill\qed

\par
Recall from Section \ref{sec:intro}  the constant $\delta_k$ for non-negative integers $k$, and that $T(X)$ denotes the set of positive square-free integers $D< X$.
\begin{theorem}\label{thm:MainTheorem1}
Let $A$ be a  positive square-free  integer such that $A\equiv 1$ or $25 \mod 36$. Then, 
for a non-negative integer $k$, 
\begin{equation}\label{eq:thm1:1}
\lim\inf_X\ \frac{ \Sharp\ \set{ D\in T(X) :  \rank E_D(\ratn)\le 2k } }{  \Sharp\ T(X) } \ge 
         \frac{\delta_k}{8} \cdot \prod_{p \mid A} \frac{p}{(p-1)(p+1)}.
\end{equation}
In particular, 
\begin{equation}\label{eq:thm1:2}
\lim\inf_X\ \frac{ \Sharp\ \set{ D\in T(X) :  E_D(\ratn)=\set{O} } }{  \Sharp\ T(X) } \ge 
         \frac{1}{16} \cdot \prod_{p \mid A} \frac{p}{(p-1)(p+1)}.
\end{equation}         
\end{theorem}
\proof
Let  $D\in\Ntwop(X/4A,1,12A)$. Then, $D$ is a square-free integer coprime to $A$ such that $D\equiv 1 \mod 12$.  By Lemma \ref{lem:condition}, $-AD^3$ satisfies Conditions (\ref{prop:cond1}) and (\ref{prop:cond2}) in Theorem \ref{prop:stoll}. 
Recall that $E_D$ is given by $y^2=x^3-AD^3$.
Note that if $D \equiv 1 \mod 12$, then $D^3\equiv 1 \mod 9$. 
Since $-AD^3 \equiv -A \equiv 2$ or $8 \mod 9$, by Theorem \ref{prop:stoll}, 
$$\rank E_D(\ratn) \le \dimp3 \Sel{\lambda}(E_D,K)=2 \dimp3 \Cl\big(\ratn(\sqrt{-AD^3})\big)[3]=2\,\log_3 h_3(-4AD).$$
\par
Let $m:=48A^2-4A$, and note that there is a one-to-one correspondence between $\Ntwop(X/4A,1,12A)$ and $\Ntwo(X,m,48A^2)$ given by $D \mapsto -4AD$. Then it follows that for a non-negative integer $k$, 
\begin{align}
\set{ \Delta \in \Ntwo(X,m,48A^2) &:  h_3(\Delta) \le 3^k }\notag\\
                 &\injects \set{D : -4AD \in \Ntwo(X,m,48A^2),\ \rank E_D(\ratn)\le 2k }.\label{eq:inclusion}
\end{align}
Let $h:=h_3$, and $B:=2$. Then, by Lemma \ref{lem:liminf} and Theorem \ref{thm:davenport}, given $\ep>0$, 
\begin{equation}\label{eq:liminf}
\frac{1}{\Sharp\, \Ntwo(X,m,48A^2)}\ \Sharp\set{ \Delta \in \Ntwo(X,m,48A^2) :  h_3(\Delta) \le 3^k } \ge \delta_k -\ep
\end{equation}
for all sufficiently large $X$.
Note that $\set{ D : -4AD \in \Ntwo(X,m,48A^2) }$ is contained in $T(X/4A)$.  Then, it follows that  given $\ep >0$, for all sufficiently large $X$,
\begin{align}
\frac{1}{\Sharp\ T(X/4A)}\ &\Sharp\ \set{ D \in T(X/4A) : \rank E_D(\ratn)\le 2k}\notag\\
&\ge \frac{1}{\Sharp\ T(X/4A)}\ 
        \Sharp\ \set{ D : -4AD\in\Ntwo(X,m,48A^2),\ \rank E_D(\ratn)\le 2k }\notag\\
        &\ge \frac{1}{\Sharp\ T(X/4A)}\ 
        \Sharp\ \set{ \Delta \in\Ntwo(X,m,48A^2):\ h_3(\Delta)\le 3^k }\quad\text{by (\ref{eq:inclusion})} \notag\\        
                 & \ge \frac{\Sharp\ \Ntwo(X,m,48A^2)}{\Sharp\ T(X/4A)}\cdot (\delta_k-\ep)
                        \quad\text{ by (\ref{eq:liminf}).}
                 \label{eq:the-lower-bound}
                 \end{align}
By \cite[Proposition 2]{nakagawa:1988}, we find 
\begin{equation}\label{eq:theconstant}
\lim_{X\to\infty}\ %
         \frac{\Sharp\ \Ntwo(X,m,48A^2)}{\Sharp\ T(X/4A)}
                =\frac{1}{8}\ \prod_{p \mid A} \frac{p}{(p-1)(p+1)},
                \end{equation}  
and this proves  (\ref{eq:thm1:1}).
\par
Let $E'/\ratn$ be an elliptic curve given by $y^2=x^3+B$ such that $B$ is an integer not equal to $-432$, $1$, a cube, or a square.  Then, it is well-known that the torsion subgroup of $E'(\ratn)$ is trivial and, hence, for all  but finitely many square-free integers $D$, the torsion subgroup of $E_D(\ratn)$ is trivial.  Therefore, (\ref{eq:thm1:2}) follows from (\ref{eq:thm1:1}) with $k=0$.
\hfill\qed

\section{Proof of Theorem \ref{thm:MainTheorem2}}\label{sec:proof2}

Let $A$ be a  square-free integer such that $A\equiv 1$ or $25 \mod 36$.   Let $m:=48 A^2 -4A$ if $A>0$, and $m:=-4A$ if $A<0$.
Note that $A\equiv 1$ or $7 \mod 9$, and that $-AD^3 \equiv -A \equiv 2$ or $8 \mod 9$ for $D \in \Ntwop(X/4\abs{A},1,12\abs{A})$ since $D \equiv 1 \mod 12$ implies $D^3 \equiv 1 \mod 9$.  
Recall that $E_D$ is given by $y^2=x^3-AD^3$. 
By Lemma \ref{lem:condition}, if $D \in \Ntwop(X/4\abs{A},1,12\abs{A})$, then $E_D$ satisfies Conditions (a) and (b) in Theorem \ref{prop:stoll} and, hence, \\ %linebreak
$$\dimp3 \Sel{\lambda}(E_D,K) = 
        \begin{cases}
                2\ \logt h_3(-4AD) & \text{if } A>0;\\
                1+2\ \logt h_3(-4AD) & \text{if } A<0.
                \end{cases}$$  
\par
If $A>0$, then there is a one-to-one correspondence between  $\Ntwop(X/4A,1,12A)$ and $\Ntwo(X,m,48 A^2)$ given by $D \mapsto -4AD$.
If $A<0$, then there is a one-to-one correspondence between   $\Ntwop(X/4\abs{A},1,12\abs{A})$ and
$\Ntwop(X,m,48 A^2)$ given by $D \mapsto -4AD$.
Note that if $n$ is a positive integer which is a power of $3$, then $\logt n \le \frac12 (n-1)$.
Then, 
it follows that if $A > 0$, then 
\begin{align*}
\frac{ \sum_{D \in \Ntwop(X/4A,1,12A)} \dimp3 \Sel{\lambda}(E_D,K) }{ \Sharp\ \Ntwop(X/4A,1,12A) }
        &=\frac{ \sum_{-4AD \in \Ntwo(X,m,48 A^2)} \dimp3 \Sel{\lambda}(E_D,K) }{ \Sharp\ \Ntwo(X,m,48 A^2) }\\
        &=\frac{ \sum_{\Delta \in \Ntwo(X,m,48 A^2)} 2\ \logt h_3(\Delta) }{ \Sharp\ \Ntwo(X,m,48 A^2) }\\
        &\le\frac{ \sum_{\Delta \in \Ntwo(X,m,48 A^2)} 2\ \frac12(h_3(\Delta)-1) }{ \Sharp\ \Ntwo(X,m,48 A^2) }
         \to 1 \\ &\qquad\text{ as $X \to \infty$, by Theorem \ref{thm:davenport}}.
\end{align*}    
If $A<0$, then 
\begin{align*}
\frac{ \sum_{D \in \Ntwop(X/4\abs{A},1,12\abs{A})} \dimp3 \Sel{\lambda}(E_D,K) }{ \Sharp\ \Ntwop(X/4\abs{A},1,12\abs{A}) }
        &\le\frac{ \sum_{\Delta \in \Ntwop(X,m,48 A^2)} 1+2\,\frac12(h_3(\Delta)-1)}{ \Sharp\ \Ntwop(X,m,48 A^2) }\\
        &\to  \frac43\quad\text{ as $X\to\infty$}.
        \end{align*}
Since the $\lambda$-Selmer rank over $K$ bounds from above the Mordell-Weil rank over $\ratn$, we have proved 
\begin{theorem}\label{thm:MainTheorem2}
Let $E/\ratn$ be an elliptic curve given by $y^2=x^3-A$ where $A$ is a square-free integer such that 
 $A \equiv 1$ or $25 \mod 36$.
\par
If $A>0$, then 
$$
\underset{X\to \infty}{\lim\sup}\ \frac{ \sum_{D \in \Ntwop(X,1,12A)} \rank( E_D(\ratn)) }{ \Sharp\ \Ntwop(X,1,12A) } \le 1.$$
\par
If $A<0$, then 
$$
\underset{X\to \infty}{\lim\sup}\  \frac{ \sum_{D \in \Ntwop(X,1,12\abs{A})} \rank( E_D(\ratn)) }{ \Sharp\ \Ntwop(X,1,12\abs{A}) } \le \frac43.$$
\end{theorem}

%\paragraph{}
%\begin{quote}
%\textbf{Acknowledgments:}\quad I thank Professor Michael Stoll for reading this preprint.
%\end{quote}

 \bibliography{mordell}

\begin{thebibliography}{10}

\bibitem{chang-thesis:2004}
S.~{}Chang.
\newblock { \rm On the arithmetic of twists of superelliptic curves, \it
  preprint\rm\ (2005)}.

\bibitem{frey:1985}
G.~Frey.
\newblock {\rm A relation between the value of $L$-series of the curve :
  $y^2=x^3-k^3$ in $s=1$ and its Selmer group, Arch. Math.\bf\ 45 \rm\ (1985),
  232--238}.

\bibitem{goldfeld:1979}
D.~Goldfeld.
\newblock {\em \rm Conjectures on elliptic curves over quadratic fields, in
  Number theory, Carbondale 1979 (Proc.~Southern Illinois Conf., Southern
  Illinois Univ., Carbondale, Ill., 1979), Lecture Notes in Math. \bf 751\rm,
  108--118, Springer-Verlag, Berlin}.

\bibitem{heath:2004}
D.R. Heath-Brown.
\newblock {\rm The average analytic rank of elliptic curves, Duke Math J. \bf
  122\ \rm (2004), no. 3, 591–623}.

\bibitem{heath-brown:1994}
D.R. Heath-Brown.
\newblock {\rm The size of Selmer groups for the congruent number problem II,
  Invent. Math. \bf 118\ \rm (1994), 331--370}.

\bibitem{lieman:1994}
D.~\rm {}Lieman.
\newblock { \rm Nonvanishing of $L$-series associated to cubic twists of
  elliptic curves, \rm Ann. of Math. \bf 140\rm\ (1994), 81--108}.

\bibitem{nakagawa:1988}
J.~Nakagawa and {K.~Horie}.
\newblock {Elliptic curves with no rational points, Proc. AMS.\bf\ 104 \rm\
  (1988), 20--24}.

\bibitem{silverman:1986}
J.~H. Silverman.
\newblock {\em {\rm The Arithmetic of Elliptic Curves, Springer, 1986}}.

\bibitem{stoll:1998}
M.~Stoll.
\newblock {\rm On the arithmetic of the curves $y^2=x^\ell+A$ and their
  Jacobians, J. reine angew. Math.\bf\ 501\rm\ (1998), 171--189}.

\bibitem{stoll:2002}
M.~Stoll.
\newblock {\rm On the arithmetic of the curves $y^2=x^\ell+A$, II, J. Number
  Theory \bf 93\ \rm (2002), 183--206}.

\bibitem{yu:2003}
G.~Yu.
\newblock {\rm Rank $0$ quadratic twists of a family of elliptic curves,
  Compositio Math. \bf\ 135\rm\ (2003), 331--356}.

\end{thebibliography}
 \bibliographystyle{plain}
\end{document}